\documentclass[12pt]{amsart}

\usepackage{amssymb,amsthm,amsmath}

\def\span{{\rm span}\,}
\newtheorem{thm}{\bf Theorem}
\newtheorem{lem}{\bf Lemma}
\newtheorem{cor}{\bf Corollary}
\newtheorem{pro}{\bf Proposition}
\newtheorem{defi}{\bf Definition}
\theoremstyle{remark}

\theoremstyle{definition}

\begin{document}
\title{Kaczmarz algorithm in  Hilbert space and tight frames}
\author{Ryszard Szwarc}
\date{}

\address{R.~Szwarc, Institute of Mathematics, University of Wroc\l aw, pl.\ Grunwaldzki 2/4, 50-384 Wroc\l aw, Poland}
\email{szwarc@math.uni.wroc.pl}
\thanks{Supported by European Commission Marie Curie Host
Fellowship for the Transfer of Knowledge ``Harmonic Analysis, Nonlinear
Analysis and Probability''  MTKD-CT-2004-013389
 and
KBN (Poland), Grant 2 P03A 028 25.}

\keywords{Kaczmarz algorithm, Hilbert space, Gram matrix, Bessel sequence, tight frame}
\subjclass[2000]{Primary 41A65}

\begin{abstract}
We prove that any   tight frame  in Hilbert space
can be obtained by the Kaczmarz algorithm. An explicit way of constructing
this correspondence is given. The uniqueness of
the correspondence is determined.
\end{abstract}
\maketitle
\section{Introduction}
Let $\{e_n\}_{n=0}^\infty$ be a linearly dense sequence of unit vectors in a Hilbert space $\mathcal H.$
 Define
{\setlength\arraycolsep{2pt}
\begin{eqnarray*}
x_0&=& \langle x,e_0\rangle e_0,\\
x_n&=& x_{n-1}+ \langle x- x_{n-1},e_n\rangle e_n.
\end{eqnarray*}
The formula is called the Kaczmarz algorithm (\cite{K}).

It can be shown   that if vectors the $g_n$ are given by the recurrence
relation
\begin{equation}\label{gn}
g_0=e_0,\quad g_n=e_n-\sum_{i=0}^{n-1}\langle e_n,e_i\rangle g_i
\end{equation}
then $g_0$ is orthogonal to $g_n,$ for any $n\ge 1$ and
\begin{equation}\label{xn}x_n=\sum_{i=0}^n\langle x,g_i\rangle e_i.
\end{equation}

By (\ref{gn}) the  vectors $\{g_n\}_{n=0}^\infty$ are linearly dense in
$\mathcal H.$
Also by definition of the algorithm
  the vectors $x-x_n$ and $e_n$ are orthogonal to each
other. Hence
\begin{eqnarray}
\|x\|^2&=&\|x-x_0\|^2+|\langle x, g_0\rangle|^2,\nonumber\\
\|x-x_{n-1}\|^2&=&\|x-x_n\|^2+ |\langle x, g_n\rangle|^2, \quad n\ge
1.\label{tight}
\end{eqnarray}

For $n\ge 1 $ let $S_n$ denote the finite dimensional operator defined by the rule
\begin{equation}
S_ny=\sum_{j=0}^{n}\langle y,e_j\rangle g_j,\quad y\in \mathcal H.
\end{equation}
Observe that the formulas (\ref{gn}) and (\ref{xn}) can be restated as
\begin{eqnarray}
(I-S_{n-1})e_n&=&g_n\\
(I-S_n^*)x&=&x-x_{n}.\label{impor}
\end{eqnarray}
Moreover by (\ref{tight}) it follows that
\begin{equation}\label{finite}
\|x-x_n\|^2=\|(I-S_n^*)x\|^2=\|x\|^2-\sum_{j=0}^{n}|\langle x,g_j\rangle|^2.
\end{equation}
In particular
\begin{equation}\label{subframe}
\sum_{n=0}^{\infty}|\langle x,g_n\rangle|^2\le \|x\|^2, \quad x\in \mathcal
H.
\end{equation}

The sequence $\{e_n\}_{n=0}^\infty$ is called effective if $x_n\to x$ for
any $x\in \mathcal H.$ By virtue of (\ref{finite}) this is equivalent to
$\|x\|^2=\sum_{n=0}^\infty ||\langle x,g_n\rangle|^2$ for any $x\in \mathcal H,$
which means $\{g_n\}_{n=0}^\infty$ is a normalized tight frame.
We refer to \cite{KM} and \cite{HS} for more information on Kaczmarz algorithm.

{\bf Acknowledgement.} I thank Wojtek Czaja and Pascu Gavruta for pointing my attention
to Lemma 3.5.1 of \cite{c}.
\section{Bessel sequences}
\begin{defi} A sequence of vectors $\{g_n\}_{n=0}^\infty$ in a Hilbert space $\mathcal H$
 will be called a
 Bessel sequence if  (\ref{subframe}) holds. The sequence $\{g_n\}_{n=0}^\infty$
 will be called a normalized Bessel sequence if in addition $\|g_0\|=1.$
 \end{defi}
Let $P_n$ denote the orthogonal projection onto $e_n^\perp$ the orthogonal complement
to the vector $e_n.$ By \cite[(1)]{HS} we have
\begin{eqnarray}
I-S_n^*&=&P_nP_{n-1}\ldots P_0,\label{ker0}\\
I-S_n&=&P_0\ldots P_{n-1}P_n.\label{ker}
\end{eqnarray}
\begin{thm}
For any normalized Bessel sequence $\{g_n\}_{n=0}^\infty$ in a Hilbert space $\mathcal H$
there exists a   sequence $\{e_n\}_{n=0}^\infty$ of unit vectors such that
(\ref{gn}) holds. In other words any normalized Bessel sequence can be obtained
through Kaczmarz algorithm.
\end{thm}
\begin{proof}
We will construct the sequence $\{e_n\}_{n=0}^\infty$ recursively. Set $e_0=g_0.$ Assume
the unit vectors $e_1,\ldots,e_{N-1}$ have been constructed such that the
formula (\ref{gn}) holds for $n=0,\ldots, N-1.$
We want to solve in $y$   the equations
\begin{equation}\label{en}(I-S_{N-1})y=g_N,\quad \|y\|=1.\end{equation}
By (\ref{ker}) we have $(I-S_{N-1})e_{N-1}=0,$ i.e. the operator
$I-S_{N-1}$ admits nontrivial kernel. Hence the solvability of
(\ref{en}) is equivalent to that of
\begin{equation}\label{enl}(I-S_{N-1})y=g_N,\quad \|y\|\le 1.
\end{equation}
By the Fredholm alternative the equation $(I-S_{N-1})y=g_N$ is solvable if
and only if $g_N$ is orthogonal to $\ker (I-S_{N-1}^*).$ We will check
that this
 condition holds. Let $x\in \ker (I-S_{N-1}^*).$ Then by (\ref{finite}) and
(\ref{subframe})
we have
$$0=\|(I-S_{N-1}^*)x\|^2=\|x\|^2-\sum_{j=0}^{N-1}|\langle
x,g_j\rangle|^2
\ge \sum_{j=N}^\infty\langle
x,g_j\rangle|^2.
$$
In particular $\langle x,g_N\rangle =0,$ i.e. $g_N\perp \ker (I-S_{N-1}^*).$

Let $y$ denote the unique solution to $$(I-S_{N-1})y=g_N,\qquad y\perp \ker
(I-S_{N-1}).$$ The proof will be complete if we show
that $\|y\|\le 1.$ Again by the Fredholm alternative we have $y\in {\rm Im}\,
(I-S_{N-1}^*).$ Let $y=(I-S_{N-1}^*)x$ for some $x\in \mathcal H.$ We may
assume that $x\perp \ker(I-S_{N-1}^*).$ In particular $\langle
x,g_0\rangle =0,$ as   (\ref{ker0}) yields  $g_0\in \ker (I-S_{N-1}^*).$
By (\ref{finite}) we have
$$\|y\|^2=\|(I-S_{N-1}^*)x\|^2=\|x\|^2-\sum_{j=1}^{N-1}|\langle
x,g_j\rangle|^2.$$
One the other hand
$$\|y\|^2=\langle x,(I-S_{N-1})y\rangle =\langle x, g_N\rangle.$$
Therefore
$$\|y\|^2-\|y\|^4=\|x\|^2-\sum_{j=1}^{N}|\langle
x,g_j\rangle|^2\ge 0,$$
which implies $\|y\|\le 1.$
\end{proof}
\begin{cor} For any normalized tight frame $\{g_n\}_{n=0}^\infty$ in a Hilbert space $\mathcal H$
there exists an effective   sequence $\{e_n\}_{n=0}^\infty$ of unit vectors such that
(\ref{gn}) holds, i.e.  any normalized tight frame can be obtained
through Kaczmarz algorithm.
\end{cor}

 For a sequence $\{e_n\}_{n=0}^\infty$ of unit vectors the normalized
Bessel sequence $\{g_n\}_{n=0}^\infty$ is
determined uniquely. However a given normalized
Bessel sequence may correspond to many sequences of unit vectors due to
two reasons. First of all for  certain $N$ the dimension of the space $\ker(I-S_{N-1})$ may
exceed 1.
 Secondly, if we fix a unit vector $u$ in $\ker(I-S_{N-1})$
the vector $e_N$ can be defined as $e_N=y+\lambda u$ for any complex
number such that $|\lambda|^2+\|y\|^2=1.$
In what follows we will indicate properties which guarantee one to one
correspondence between $\{e_n\}_{n=0}^\infty$ and $\{g_n\}_{n=0}^\infty.$

\begin{defi} A sequence of unit vectors $\{e_n\}_{n=0}^\infty$ will be called
stable if the vectors $ \{e_n\}_{n= N}^\infty $ are linearly dense
for any $N.$
A normalized Bessel sequence $\{g_n\}_{n=0}^\infty$ will be called stable
if the vectors $\{g_0\}\cup\{g_n\}_{n= N}^\infty $ are linearly dense
for any $N.$
\end{defi}
\begin{pro} Let  sequences $\{e_n\}_{n=0}^\infty$ and
$\{g_n\}_{n=0}^\infty$ satisfy (\ref{gn}). The sequence
$\{g_n\}_{n=0}^\infty$ is stable if and only if $\{e_n\}_{n=0}^\infty$ is
stable and $\langle e_n,e_{n+1}\rangle \neq 0 $  for any $n\ge 0.$
\end{pro}
\begin{proof}
Assume $\{g_n\}_{n=0}^\infty$ is stable. First we will show
that the kernel of $I-S_{N-1}$ is one
dimensional and thus consists of the multiples of the vector $e_{N-1}$
(see (\ref{ker}).  Assume for a contradiction that
$\dim\ker (I-S_{N-1})\ge 2.$ By the Fredholm alternative we get
$\dim\ker (I-S_{N-1}^*)\ge 2.$ Hence there exists a nonzero vector $x$
such that $x\perp g_0$ and $(I-S_{N-1}^*)x=0.$ By (\ref{tight}) we obtain
$$\|x\|^2=\sum_{n=1}^{N-1}|\langle x,g_n\rangle |^2.$$
This and the  condition (\ref{subframe}) imply that $x$ is orthogonal to all the vectors $g_0$
and
$\{g_n\}_{n=N}^\infty,$ which contradicts the stability assumption.

Assume $\langle e_{N-1},e_{N}\rangle =0$ for some $N\ge 1.$ Then by
(\ref{ker}) we have
$e_{N-1},\,e_N\in \ker (I-S_{N-1})$ which is a contradiction as the kernel is
one dimensional.

Concerning stability of $\{e_n\}_{n=0}^\infty$ assume a vector $y$ is
orthogonal to all the vectors
$\{e_n\}_{n=N}^\infty.$ In particular $y$ is orthogonal to $e_N.$ Since
$\ker(I-S_N)=\mathbb{C}e_N$ by the Fredholm alternative $y$ belongs to
${\rm Im}\,(I-S_N^*).$ Let $y=(I-S_N^*)x$ for some $x\in \mathcal H.$
We may assume that $x\perp g_0$ as $g_0\in \ker (I-S_N^*).$
By (\ref{ker0}), since $y$ is orthogonal to $e_n$ for $n\ge N,$ we get
$y=(I-S_n^*)x=(I-S_N^*)x$ for $n\ge N.$ On the other hand by (\ref{xn}) and
(\ref{impor}) we
obtain that $\langle x,g_n\rangle =0$ for $n>N+1.$ Since $x\perp g_0$ by stability
assumptions we obtain $x=0$ and thus $y=0.$

For the converse implication assume $\{e_n\}_{n=0}^\infty$ is stable
and $\langle e_n,e_{n+1}\rangle \neq 0.$
By  the inequality (see \cite{KM})
$$\|x-x_n\|\ge |\langle e_{n-1},e_n\rangle| \|x-x_{n-1}\|$$ we get
that $x-x_n\neq 0$ for any $x\perp e_0.$ Since $x-x_n=(I-S_n^*)x$ the kernel of $I-S_n^*$
consists only of the multiples of $e_0=g_0.$

Let $x$ be orthogonal to $\{g_0\}\cup\{g_n\}_{n\ge N+1}$ for some $N\ge 1.$
By (\ref{xn}) we obtain that $x_n=x_N$ for $n\ge N.$ By the definition of
the Kaczmarz algorithm we get   $x-x_N\perp e_n$ for $n\ge N+1.$
Now stability of $\{e_n\}_{n=0}^\infty$ implies $x-x_N=0.$ By
  (\ref{impor}) we obtain $(I-S_N^*)x=0.$ This implies $x=0$ since
  the kernel is one dimensional and consists of the multiples of $g_0.$
\end{proof}
For sequences $\{e_n\}_{n=0}^\infty $ and $\{\sigma_ne_n\}_{n=0}^\infty , $
where $\sigma_n$ are complex numbers of absolute value 1, the Kaczmarz
algorithm coincide. Therefore we will restrict our attention to {\em admissible} sequences
of unit vectors $\{e_n\}_{n=0}^\infty$ such that $\langle
e_n,e_{n+1}\rangle \ge 0.$
\begin{thm}
Let $\{g_n\}_{n=0}^\infty$ be a stable normalized Bessel sequence. Then there
exists a unique admissible sequence $\{e_n\}_{n=0}^\infty $ of unit
vectors such that (\ref{gn}) holds. Moreover the sequence  $\{e_n\}_{n=0}^\infty $
is stable.
\end{thm}
\begin{proof}
The proof will go by induction. The vector $e_0$ is determined   by
$e_0=g_0.$ Assume the vectors $e_0,\ldots, e_{N-1}$ were determined
uniquely. We have to show that the problem
$$
(I-S_{N-1})y=g_N,\quad \|y\|=1,\quad \langle y,e_{N-1}\rangle \ge 0$$
has the unique solution $y.$

By the proof of Proposition 2 the kernel of $I-S_{N-1}$ is one
dimensional and thus consists of the multiples of the vector $e_{N-1}.$
 By the proof of Theorem 1 there exists
 the unique solution $y_N$ to the problem
$$(I-S_{N-1})y=g_N,\quad y\perp \ker (I-S_{N-1})$$
 and  $\|y_N\|\le 1.$ Moreover by
this proof
$\|y_N\|=1$ if and only if $$\|x\|^2-\sum_{j=1}^{N}|\langle
x,g_j\rangle|^2=0,$$ where $y_N=(I-S_{N-1}^*)x$ and $x\perp \ker(I-S_{N-1}^*).$
This leads to a contradiction because by inequality (\ref{subframe})  we get that
$x$ is orthogonal to all the vectors $g_0$
and
$\{g_n\}_{n=N}^\infty.$ Hence $\|y_N\|<1.$

At this stage we know that any solution to the equation
$$(I-S_{N-1})y=0$$ is of the form
$$y=y_N+\lambda e_{N-1},\quad \lambda\in \mathbb{C}$$
because $\ker (I-S_{N-1})=\mathbb{C}e_{N-1}.$
Since $\|y_N\|<1$ and $y_N\perp e_{N-1}$ there exists a unique solution
$y$
satisfying $\|y\|=1$ and $\langle y, e_{N-1}\rangle \ge 0$ namely
the one corresponding to $\lambda=\sqrt{1-\|y_N\|^2}.$
\end{proof}
\begin{cor} Let $\{g_n\}_{n=0}^\infty$ be a stable normalized tight frame. Then there
exists a unique admissible effective sequence $\{e_n\}_{n=0}^\infty $ of unit
vectors such that (\ref{gn}) holds. Moreover the sequence  $\{e_n\}_{n=0}^\infty $
is stable.
\end{cor}
\section{Algorithm}
The proof of Theorem 1 can also be given by using Gram matrix of the
sequence $\{g\}_{n=0}^\infty.$ This argument can be used for
constructing an underlying sequence of unit vectors
$\{e_n\}_{n=0}^\infty.$ This will be done below.

The following lemma follows from Lemma 1 it follows from  \cite[Lemma 3.5.1]{c}.

\begin{lem}\label{iff} The collection $\{g_n\}_{n=0}^\infty$ is a   Bessel sequence
if and only if the Gram matrix $G=\{\langle g_i,g_j\rangle\}_{i,j=0}^\infty $ corresponds to a
contraction operator on $\ell^2(\mathbb{N}).$ The sequence $\{g_n\}_{n=0}^\infty$
is a tight frame if and only if it is linearly dense and  $G$ corresponds to a projection
on $\ell^2(\mathbb{N}).$
\end{lem}

Let $\{e_n\}_{n=0}^\infty$ be a sequence of unit vectors in a Hilbert
space $\mathcal H$ and let $\{g_n\}_{n=0}^\infty$ be the corresponding
normalized Bessel sequence.
Let $M$   be the strictly lower triangular part of the Gram  matrix of
the sequence $\{e_n\}_{n=0}^\infty$  and
  $U$   strictly lower triangular matrix defined by $$(I+U)(I+M)=I.$$ By \cite{HS}
the matrix $U$ is a contraction on the Hilbert space $\ell^2(\mathbb{N}).$
\begin{lem}\label{if}
For any $i,j$ we have
$$\langle g_i,g_j\rangle =\langle
(I-UU^*)\delta_j,\delta_i\rangle_{\ell^2(\mathbb{N})}$$
\end{lem}
\begin{proof}
Let
$$M=\begin{pmatrix}
0 & 0 & 0&0 & 0&\ldots \\
m_{10} & 0 &0& 0&0 & \ldots \\
m_{20} & m_{21} & 0&0 &0&\ldots\\
m_{30}&m_{31}&m_{32}& 0&0&\ldots\\
\vdots&\vdots & \vdots &\vdots&\ddots&\ddots
\end{pmatrix},\quad U=\begin{pmatrix}
0 & 0 & 0&0 & 0&\ldots \\
c_{10} & 0 &0& 0&0 & \ldots \\
c_{20} & c_{21} & 0&0 &0&\ldots\\
c_{30}&c_{31}&c_{32}& 0&0&\ldots\\
\vdots&\vdots & \vdots &\vdots&\ddots&\ddots
\end{pmatrix}.$$
 By \cite[(6)]{HS} we have
\begin{equation}\label{equiv}
g_i=e_i+\sum_{k=0}^{i-1}c_{ik}e_k.
\end{equation}
Set $c_{nn}=m_{nn}=1$ and let $1\le k\le n.$ By taking inner product with $g_j$
in (\ref{equiv}) we get
\begin{multline*}
\langle g_i,g_j\rangle =
\sum_{k=0}^ic_{ik}\sum_{l=0}^j\overline{c}_{jl}\langle e_k,e_l\rangle
\\=\langle (I+U)(I+M+M^*)(I+U^*)\delta_j,\delta_i\rangle
_{\ell^2(\mathbb{N})}.
  \end{multline*}
Taking into account relations between the matrices $M$ and $U$ gives
\begin{equation}\label{MU}(I+U)(I+M+M^*)(I+U^*)=I-UU^*.
\end{equation}
The product of these matrices is well defined since $U^*$ leaves the space
$F(\mathbb{N})={\rm span}\, \{\delta_n\mid n\ge 0\}
)$ invariant. Therefore
$$\langle g_i,g_j\rangle=\langle
(I-UU^*)\delta_j,\delta_i\rangle_{\ell^2(\mathbb{N})}.$$
\end{proof}
{\bf Remark.} Lemma \ref{if} can be used to give a shorter and simpler proof of
Theorem 1 in \cite{HS}. Indeed, by Lemma \ref{iff}, a linearly
dense sequence $\{g_n\}_{n=0}^\infty$ constitutes a normalized tight frame
 if and only if the matrix $G=\{\langle g_i,g_j\rangle\}_{i,j=0}^\infty $ is a
projection. In view of Lemma \ref{if} the latter is equivalent to $U$ being a
partial isometry. Moreover in this case we have
$$\dim{\mathcal H}=\sum_{n=0}^\infty \|g_n\|^2={\rm Tr}\, (I-UU^*).$$

We are ready now to give an alternative proof of Theorem 1. Let
$\{g_i\}_{i=0}^\infty$ be a normalized Bessel sequence. By Lemma \ref{iff}
the matrix $A:=I-G$ is positive definite.
 Moreover
$A(0,i)=A(i,0)=0$ because $\|g_0\|=1$ and $g_0\perp g_i$ for $i\ge 1.$
Let $\tilde{A}$ denote the truncated matrix obtained from $A$ by removing
the first row and the first column. Clearly $\tilde{A}$ corresponds to a
positive definite contraction on $\ell^2(\mathbb{N}^+).$
The next lemma is  probably known and provides infinite dimensional version
of the so called Cholesky
decomposition of positive definite matrices.
\begin{lem}\label{three} For any positive definite matrix $B=\{b(i,j)\}_{i,j=1}^\infty$ there exists a
lower triangular matrix $V=\{v(i,j)\}_{i,j=1}^\infty$ such that $B=VV^*.$
\end{lem}
\begin{proof}
By the well known fact there exists a Hilbert space $\mathcal M$ and a
linearly dense
sequence of vectors $\{h_i\}_{i=1}^\infty$ in $\mathcal M$ such that
$$b(i,j)=\langle h_i,h_j\rangle.$$
By applying the Gram-Schmid procedure to this sequence we obtain an
orthonormal sequence $\{\eta_i\}_{i=1}^N,$ where $N=\dim \mathcal M,$
such
that $h_i\in {\rm span\,}\{\eta_1,\eta_2,\ldots, \eta_i\}$ for $i< N+1.$
In particular there are coefficients $v_{ik}$ for $i\ge k$ and $N+1> k$ for which
we have
$$h_i=\sum_{k=1}^iv_{ik}\eta_k.$$
Set $v_{ik}=0$ for $i>k$ and for $k>N.$ Then
$$
b(i,j)=\langle h_i,h_j\rangle =
\sum_{k,l=1}^{\min(i,j)}v_{ik}\overline{v}_{jk}=
(VV^*)(i,j).
$$
\end{proof}
By Lemma \ref{three} there is a lower triangular matrix $V=\{v_{ij}\}_{i,j=1}^\infty$ such that
$\widetilde{A}=VV^*.$ Let $U=\{c_{ij}\}_{i,j=0}^\infty$
 be the strictly lower triangular matrix obtained from $V$ by adding a
 zero row and a zero column, i.e.
 $$c_{ij}=\begin{cases}v_{i-1,j-1}& ij>0,\\
 0 &ij=0.
 \end{cases}$$
 In this way we obtain
 \begin{equation}\label{I-G}I-G=A=UU^*.
 \end{equation}
 Since $G$ corresponds to a contraction on $\ell^2(\mathbb{N})$ so does
 $U.$
 Let $M=\{m_{ij}\}_{i,j=0}^\infty$ be the strictly lower triangular matrix
 determined by $(I+M)(I+U)=I.$ Set $m_{ii}=1$ and define (cf. (\ref{gn}))
 $$e_i=\sum_{k=0}^{i}m_{ik}g_k.$$
 We claim that $e_i$ are unit vectors and moreover $\langle e_i,e_j\rangle
 =m_{i,j}$ for $i\ge j.$ This will give (\ref{gn}) and thus conclude another
 proof of Theorem 1.
By (\ref{I-G}) we have
\begin{multline*}
\langle e_i,e_j\rangle =
\sum_{k=0}^i\sum_{l=0}^jm_{ik}\overline{m}_{jl}\langle g_k,g_l\rangle
=\sum_{k=0}^i\sum_{l=0}^jm_{ik}\overline{m}_{jl}(I-UU^*)(k,l)\\=
(I+M)(I-UU^*)(I+M^*)(i,j)
\end{multline*}
On the other hand (\ref{MU}) yields $$(I+M)(I-UU^*)(I+M^*)=I+M+M^*.$$
In particular for $i\ge j$ we obtain
$$\langle e_i,e_j\rangle
 =\begin{cases}m_{i,j}& i>j,\\
 1&i=j.\end{cases}
 $$

This way of proving Theorem 1 provides an algorithm for constructing a
sequence of unit vectors $\{e_n\}_{n=0}^\infty$ for a given normalized
Bessel sequence $\{g_n\}_{n=0}^\infty.$ Indeed, it suffices to determine
an algorithm for proving Lemma \ref{three}. When $B$ is strictly positive definite then
the solution can be given by the so called Cholesky algorithm. When $B$ is not necessarily
strictly positive definite
this algorithm fails and w ehave to find a different way of constructing the decomposition.

We will construct a sequence of indices $\{n_k\}_{k=1}^N$ in the following way.
Let $n_1$ be the smallest number $i$ such that $b_{ii}>0.$
If such number does not exist then $B=0.$
Assume $n_1,n_2,\ldots, n_k$ have been constructed in such a way
that the determinant
$$\Delta_k=\det (b_{n_in_j})_{i,j=1}^k> 0.$$
Then let $n_{k+1}$ be the smallest number such that
$$\det (b_{n_in_j})_{i,j=1}^{k+1}> 0.$$
If such number does not exist the procedure terminates and $N=n_k.$

The matrix $B$ gives rise to a positive definite
hermitian form on the space $F(\mathbb{N}_+)=\span\{\delta_n\,|\,n\ge 1\}$ by the rule
$$\langle x,y\rangle =\sum_{i,j=1}^\infty b(i,j)x_i\overline{y}_j.$$
\begin{lem}
For any $n$ there exist $i$ with
$n_{i}\le n<n_{i+1}$ and numbers $\lambda_{nk}$ for $1\le k\le i,$ such that
\begin{equation}\label{equ}
\left \langle \delta_n-\sum_{k=1}^i\lambda_{nk}\delta_{n_k},\delta_m\right
\rangle =0, \quad m\ge 1.
\end{equation}
\end{lem}
\begin{proof}
If $n=n_i$ for some $i,$ then statement follows as
$\delta_n=\delta_{n_i}.$ Otherwise we have $n_i<n<n_{i+1}$ for some
$i.$ By plugging in $m=n_1,n_2,\ldots, n_i$ to (\ref{equ}) we obtain
a system of linear equations
$$\sum_{k=1}^i\lambda_{nk}b_{n_kn_l}=b_{nn_l},\quad l=1,2,\ldots, i.$$
The main determinant of this system is $\Delta_i.$ Therefore the system
has the unique solution $\lambda_{n,1},\ldots,\lambda_{n,i}.$
By definition of $n_{i+1}$ we have
$$ \begin{vmatrix}
b_{n_1n_1}&b_{n_1n_2}&\cdots & b_{n_1n_i}&b_{n_1n}\\
b_{n_2n_1}&b_{n_2n_2}&\cdots & b_{n_2n_i}&b_{n_2n}\\
\ \vdots &\ \vdots &\cdots &\ \vdots&\vdots\\
b_{n_in_1} &b_{n_in_2}  & \cdots & b_{n_in_i}& b_{n_in}\\
b_{nn_1}& b_{nn_2} & \cdots & b_{nn_i}& b_{nn}
\end{vmatrix}
=0.
$$
As $\Delta_i>0$ the first $i$ rows of this matrix are linearly independent.
Therefore the last row of this matrix is a linear combination of the first
$i.$ The coefficients must coincide with $\lambda_{n1},\ldots,\lambda_{ni}.$
In particular considering the last entry of the rows gives
$$\sum_{k=1}^i\lambda_{nk}b_{n_kn}=b_{nn}.$$
This is equivalent to (\ref{equ}) with $m=n.$

Since (\ref{equ}) is valid for $m=n,n_1,\ldots, n_i$ then
$$\left\langle \delta_n-\sum_{k=1}^i\lambda_{nk}\delta_{n_k},
\delta_n-\sum_{k=1}^i\lambda_{nk}\delta_{n_k}\right\rangle =0.$$
By Schwarz inequality this implies (\ref{equ}) for any $m.$
\end{proof}
Define the sequence of vectors $\{\eta_i\}_{i=1}^N$ by the formula
$$
\eta_1={1\over \sqrt{\Delta_1}}\delta_{n_1},\qquad\eta_i=
{1\over \sqrt{\Delta_{i-1}\Delta_i}}
\begin{vmatrix}
b_{n_1n_1}&b_{n_1n_2}&\cdots & b_{n_1n_i}\\
b_{n_2n_1}&b_{n_2n_2}&\cdots & b_{n_2n_i}\\
\ \vdots &\ \vdots &\cdots &\ \vdots\\
b_{n_{i-1}n_1}&b_{n_{i-1}n_2}&\cdots &b_{n_{i-1}n_i}\\
\delta_{n_1} &\delta_{n_2}  & \cdots & \delta_{n_i}
\end{vmatrix}
 $$
It can be checked easily that
\begin{equation}\label{orth}
\langle \eta_i,\eta_j\rangle =\delta_i^j.
\end{equation}
Obviously from the definition we have
$$\eta_{i}={\sqrt{\Delta_i}\over \sqrt{\Delta_{i-1}}}\delta_{n_i}+
\sum_{k=1}^{i-1}\alpha_{ik}\delta_{n_k}$$
for some explicitly given coefficients $\alpha_{ik}.$ Therefore
\begin{equation}\label{orthinv}
\delta_{n_i}=\sum_{k=1}^i\beta_{ik}\eta_k,
\end{equation}
for some coefficients $\beta_{ik}.$
By   (\ref{equ}) and (\ref{orthinv}) we get
that
for any $n$ there exist $i$ with
$n_{i}\le n<n_{i+1}$ and numbers $v_{nk}$ for $1\le k\le i,$ such that
$$
\left \langle \delta_n-\sum_{k=1}^iv_{nk}\eta_{k},\delta_m\right
\rangle =0,\quad m\ge 1.
$$
Setting $v_{nk}=0$ for $i<k\le n$ gives
\begin{equation}\label{equ1}
\left \langle \delta_n-\sum_{k=1}^nv_{nk}\eta_{k},\delta_m\right
\rangle =0,\quad m\ge 1.
\end{equation}
Therefore by (\ref{equ}) and (\ref{equ1}) we have
\begin{multline*}
0=\left\langle
\sum_{k=1}^nv_{n,k}\eta_k\,,\,\delta_m-\sum_{k=1}^nv_{mk}\eta_{k}\right\rangle=
\langle \delta_n,\delta_m\rangle -
\sum_{k=0}^{\min(n,m)}v_{nk}\overline{v}_{mk}\\
=b_{nm}-\sum_{k=0}^{\min(n,m)}v_{nk}\overline{v}_{mk}.
\end{multline*}
Therefore $B=VV^*$ where
$$
V=\begin{pmatrix}
v_{11} & 0 & 0&0 & 0&\ldots \\
v_{21} & v_{22} &0& 0&0 & \ldots \\
v_{31} & v_{32} & v_{33}&0 &0&\ldots\\
v_{41}&v_{42}&v_{43}& v_{44}&0&\ldots\\
\vdots&\vdots & \vdots &\vdots&\ddots&\ddots
\end{pmatrix}.
$$
By analyzing the entire construction we may conclude the coefficients $v_{nk}$
can be computed   in an algorithmic way.
\section{Equivalent sequences}
Any sequence of unit vectors $\{e_n\}_{n=0}^\infty$ leads by Kaczmarz algorithm
to a normalized Bessel sequence $\{g_n\}_{n=0}^\infty.$
\begin{defi}
Two normalized Bessel sequences $\{g_n\}_{n=0}^\infty$ and $\{g'_n\}_{n=0}^\infty$
will be called equivalent if there is a unitary operator $V$ such that
$g'_n=Vg_n$ for $n\ge 0.$
Similarly two sequences $\{e_n\}_{n=0}^\infty$ and $\{e_n'\}_{n=0}^\infty$ of unit vectors
will be called  equivalent if there is a unitary operator $V$ such that
$e'_n=Ve_n$ for $n\ge 0.$
 \end{defi}
It is easy to see that if the sequences $\{e_n\}_{n=0}^\infty$ and $\{e_n'\}_{n=0}^\infty$
are equivalent so are the corresponding sequences of normalized Bessel
sequences $\{g_n\}_{n=0}^\infty$ and $\{g_n'\}_{n=0}^\infty,$ with the
same unitary operator $V.$ The converse is not true, as the normalized
Bessel sequences do not correspond to sequences of unit vectors in
one-to-one fashion. Nevertheless
by Lemma 1,  the equivalence relation between normalized Bessel sequences can be
described in terms of Gram matrices of the corresponding sequences of
unit vectors.

Assume sequences of unit vectors $\{e_n\}_{n=0}^\infty$ and $\{e_n'\}_{n=0}^\infty$
are associated with normalized Bessel sequences $\{g_n\}_{n=0}^\infty$ and $\{e_g'\}_{n=0}^\infty,$
respectively.
Let $M$ and $M'$ be the strictly lower triangular part of the Gram  matrices of
sequences $\{e_n\}_{n=0}^\infty$ and $\{e_n'\}_{n=0}^\infty,$ respectively.
Let $U$ and $U'$ be strictly lower triangular matrices defined by $$(I+U)(I+M)=(I+U')(I+M')=I.$$ By \cite{HS}
the matrices $U$ and $U'$ are contractions on the Hilbert space $\ell^2(\mathbb{N})$

\begin{cor}
The sequences $\{g_n\}_{n=0}^\infty$ and $\{g_n'\}_{n=0}^\infty$ are
 equivalent if and only if $UU^*=U'U'^*.$
\end{cor}
\begin{proof}
By Lemma 1 we get that $UU^*=U'U'^*$ if and only if $\langle
g_i,g_j\rangle =\langle
g'_i,g'_j\rangle$ for any $i,j\ge 0.$ Obviously the latter, along with
the linear density of vectors $\{g_n\}_{n=0}^\infty $ and
$\{g'_n\}_{n=0}^\infty ,$ is equivalent
to the existence a unitary operator $U$ such that $g'_i=Ug_i.$
\end{proof}

{\bf Remark.} It would be of interest to determine when two sequences of
unit vectors $\{e_n\}_{n=0}^\infty$ and $\{e_n'\}_{n=0}^\infty$ lead
to the same normalized Bessel sequence $\{g_n\}_{n=0}^\infty .$

\end{document}